\documentclass[12pt]{amsart}
\usepackage{amsmath, amsfonts, amssymb, amscd}
\newtheorem{theo}{Theorem}[section]
\newtheorem{lem}[theo]{Lemma}
\newtheorem{prop}[theo]{Proposition}

\newtheorem{corollary}[theo]{Corollary}

\newtheorem{defi}[theo]{Definition}
\newenvironment{pf}{\noindent{\it Proof. }}{\hspace{8cm}$\square$\par\medskip}

\newenvironment{rem}{\noindent{\it Remark.}}{}
\newcommand{\C}{{\mathbb C}}
\newcommand{\R}{{\mathbb R}}
\newcommand{\Z}{{\mathbb Z}}

\renewcommand{\l}{{\mathfrak l}}
\renewcommand{\k}{{\mathfrak k}}
\newcommand{\m}{{\mathfrak m}}

\renewcommand{\a}{{\mathfrak a}}
\newcommand{\h}{{\mathfrak h}}

\newcommand{\n}{{\mathfrak n}}

\newcommand{\SU}{{\operatorname{SU}}}
\newcommand{\T}{{\operatorname{T}}}
\newcommand{\SL}{{\operatorname{SL}}}
\newcommand{\UU}{{\operatorname{U}}}

\newcommand{\SO}{{\operatorname{SO}}}
\newcommand{\Spin}{{\operatorname{Spin}}}
\newcommand{\EE}{{\operatorname{E}}}
\newcommand{\GG}{{\operatorname{G}}}
\newcommand{\Sp}{{\operatorname{Sp}}}

\newcommand{\rk}{\operatorname{rank}}

\newcommand{\rank}{\operatorname{rank}}

\newcommand{\Ad}{\operatorname{Ad}}

\title[Complex asystatic actions of compact Lie groups]
{Complex asystatic actions of compact Lie groups}
\author{Anna Gori and Fabio Podest\`{a}}
\address {Dipartimento di Matematica e Appl.\ per l'Architettura\\
Piazza Ghiberti 27\\50100 Firenze\\Italy } \email{gori@math.unifi.it}
\address{Dipartimento di Matematica e Appl.\ per l'Architettura\\
Piazza Ghiberti 27\\50100 Firenze\\Italy}
\email{podesta@math.unifi.it}
\date{\today}

\begin{document}
\maketitle
\null \vspace*{-.25in}
\begin{abstract}
In the present paper we introduce the notion of complex asystatic
Hamiltonian action on a K\"ahler manifold. In the algebraic
setting we prove that if a complex linear  group $G$ acts
complex asystatically on a K\"ahler  manifold then the $G$-orbits
are spherical. Finally we give the complete classification of
complex asystatic irreducible representations.
\end{abstract}
\section{Introduction}
\medskip
A proper action of a compact Lie group $K$ on a
manifold $M$ is called asystatic if the isotropy representation of
a principal isotropy subgroup $L$ has no non-trivial fixed vector
in the tangent space of a principal orbit $K/L$, or, equivalently,
if the quotient group  $N_K(L)/L$ is finite. Asystatic actions are
particular cases of polar actions, namely, if we endow $M$ with
any $K$-invariant Riemannian metric, there exists a properly
embedded submanifold $\Sigma$, called a section, which intersects
every $K$-orbit orthogonally; in the case of asystatic actions, a
section is provided by the fixed point set of a principal
isotropy. Asystatic actions, whose introduction goes back to
Lie (see \cite{LE}), have been treated in \cite{AA} and \cite{PT}.\\
If $K$ acts isometrically and holomorphically on a K\"ahler
manifold $M$, then the asystatic condition is rarely satisfied, at
least when the action is Hamiltonian, as the following proposition shows.
\begin{prop}\label{flag} Let $K$ acts effectively on a K\"ahler manifold $M$ by holomorphic
isometries. If the $K$-action is Hamiltonian, then it is asystatic
if and only if $K$ is semisimple and $M$ is the K\"ahler product
of a flag $K/H$ and a K\"ahler manifold which is acted on
trivially by $K$.
\end{prop}
Now, given a $K$-action on a manifold $M$ with principal isotropy type
$(L)$, the core $\overline M\subset M$ is defined as the closure
of $M^L\cap M_{\textrm{reg}}$, where $M^L$ denotes the
fixed point set of the subgroup $L$. The quotient group $N_K(L)/L$ acts
effectively on $\overline{M}$ and the inclusion
$\overline{M}\hookrightarrow M$ induces a bijection at the orbit
space level; this gives the primary motivation for referring to
the induced transformation group $((N_K(L)/L),\overline{M})$ as
the {\em reduction} of $(K,M)$ (see e.g. \cite{GS} and \cite{SS}).

 In this paper we are interested in the fixed point set of a principal
isotropy for a Hamiltonian action of a compact group of
holomorphic isometries $K$ of a K\"ahler manifold $M$; in this
case one of our results states the following
\begin{prop} \label{inequality} Let $M$ be a K\"ahler manifold endowed with an
effective Hamiltonian $K$-action. Let $(L)$ be a principal
isotropy type, $\overline M$ a core and $c$ the cohomogeneity of
the $K$-action. Then
\begin{enumerate}
\item $\dim \overline M \geq  c + \rank(K) - \rank(L)$;
\item the equality in (1) holds if and only if the group $(N_K(L)/L)^o$ is abelian;
\item if $M=V$ is a linear space and $K$ acts linearly on $V$ or $M$ is locally irreducible, then
$\rk(K) - \rk(L)\geq 1$.
\end{enumerate}
\end{prop}
It is therefore natural to investigate actions for which this
equality holds. With this notation we say that a proper $K$-Hamiltonian action on
a K\"ahler manifold $M$ by holomorphic isometries is called
{\it complex asystatic\/} (or $\C$-{\it asystatic\/}) if
$$\dim \overline M=c+\rk K-\rk L$$
or, equivalently, if $(N_K(L)/L)^o$ is abelian.\par
In the algebraic setting we prove the following theorem, which connects
the $\C$-asystatic condition with the class of actions of reductive groups with
spherical orbits; we recall here that a complex homogeneous space $N$ of a reductive
complex group $G$ is called spherical if a Borel subgroup of $G$ has an open orbit
in $N$ (see e.g.~\cite{BLV}).
\begin{theo}\label{spherical} Let $G$ be a connected complex linear algebraic reductive group acting
algebraically on a K\"ahler algebraic manifold $M.$ If the action
of a maximal compact Lie subgroup of $G$ on $M$ is Hamiltonian and
$\C$-asystatic, then every $G$-orbit is spherical.\end{theo} All
irreducible representations of a reductive algebraic group $G$
with spherical orbits have been classified by Arzhantsev
(\cite{A2}). In the last section we therefore use this result and
our Theorem~\ref{spherical} in order to determine all the
irreducible $\C$-asystatic representations. Our result is the
following
\begin{theo}\label{class} Let $K$ be a connected compact Lie group and $\rho:K\to\UU(V)$ be a complex
irreducible representation. Then $\rho$ is $\C$-asystatic if and only if $\rho(K)$ is one of the following
\begin{enumerate}
\item $\SU(n)$, $\Lambda^2(\SU(n))$, $S^2(\SU(n))$ ($n\geq 3$), $\SO(n)$, $\SU(n)\otimes \SU(m)$, $\Spin(7)$,
$\Spin(10)$, $\EE_6$;
\item $\rho(H)\cdot {\operatorname{T}}^1$, where $\rho(H)$ appears in ${\operatorname\ (1)}$, or $S^2(\UU(2))$,
${\operatorname{T}}^1\times
\Sp(n)$, ${\operatorname{T}}^1\times \GG_2$.
\end{enumerate}
\end{theo}
We note here that the actions appearing in the classification are either orbit equivalent to isotropy
representations of Hermitian symmetric spaces or are obtained from these by deleting the one-dimensional
torus.\par
The organization of the paper is as follows: definitions and basic useful facts on complex asystatic actions
are given in Section $2$, while in Section $3$ we prove our main result in the algebraic
setting. Applying our main result we give, in Section $4$, the complete classification of irreducible  complex asystatic
representations.
\section{Actions on K\"ahler manifolds and the $\C$-asystatic condition}

Throughout the following $K$ will denote a connected compact Lie
group. A proper $K$-action on a manifold $M$ is said to be {\em
asystatic} if, given a principal point $p\in M$ and its isotropy
subgroup $K_p:=L,$ the isotropy representation of $L$ on
$T_p(K\cdot p)$ has no nontrivial fixed vector.\par
If $M$ is a
K\"ahler manifold and $K$ acts on $M$ by holomorphic isometries,
the asystatic condition is rarely satisfied, at least if $K$ acts
in Hamiltonian fashion. We now give the proof of Prop\ref{flag} (for some related
see also~\cite{GP}).
\begin{proof}[Proof of Proposition \ref{flag}] Let $L$ be a principal isotropy subgroup and let
$\Sigma$ be the connected component of the fixed point set $M^L$
passing through the principal point $p$. Then $T_p\Sigma\cap T_p(K\cdot p) = \{0\}$ and
since $\Sigma$ is complex, so is $K\cdot p$. We now claim that
every $K$-orbit is complex. Indeed, $\Sigma$ intersects every $K$-orbit and for
every $q\in \Sigma$ we have that $T_q\Sigma$ is a section for the linear $K_q$-action on the normal space
$\nu_q$ of $K\cdot q$ (see e.g. \cite{PT}); moreover if $v\in \nu_q$, there exists $g\in K_q$ so that
$gv\in T_q\Sigma$, hence $Jv\in g^{-1}T_q\Sigma\subset \nu_q$ since $\Sigma$ is complex. This
means that $\nu_q$, hence $K\cdot q$ is complex.\par
Since the $K$-action
is Hamiltonian, the $Z$-action on $K\cdot p$ is also Hamiltonian,
where $Z$ denotes the connected component of the center of $K$;
this implies that $Z$ acts trivially
on $K\cdot p$, hence on $M$ by the effectiveness of the action, and therefore $K$ is semisimple.\par
If $\mu:M\to \k^*$ is a moment map for the $K$-action, where
$\k^*$ denotes the dual of the Lie algebra of $K$,
we first note that $\mu$ is constant on $\Sigma$, because $\Sigma$
is complex and therefore $\langle d\mu_q(X),Y\rangle =
\omega(X,\hat Y)|_q = 0$ for every $X\in T_q\Sigma, Y\in \k$ and $q\in \Sigma$. We put
$\mu(\Sigma) = \{\beta\}$ for some $\beta\in \k^*$. If $x\in \Sigma$, then
$K\cdot x$ is complex, hence $\mu:K\cdot x\to K\cdot\beta$ is a covering and since
$K\cdot \beta$ is simply connected, we have $K_x=K_\beta$ and every $K$-orbit in $M$ is principal.
This also implies that $\Sigma$ intersects every $K$-orbit in a single point: indeed if
$x,gx\in \Sigma$ for some $g\in K$, then $g\in K_\beta=K_x$, hence $x=gx$. The map
$\phi:K\cdot p\times \Sigma\to M$ given by $\phi(kp,x)=kx$ where $p$ is any fixed point in
$\Sigma$ is a well defined $K$-equivariant diffeomorphism, which can be proved to be an isometry
using the same arguments as in the proof of the main result in~\cite{GP}.\end{proof}

From now on we will always suppose that the $K$-action is also
Hamiltonian.\\Given a $K$-action on a manifold $M$ with principal
isotropy type $(L)$, recall that the core $\overline M\subset M$
is defined as the closure of $M^L\cap M_{\textrm{reg}}$ and it is
known (see~\cite{GS} and~\cite{SS}) that all the components of
$\overline M$ are equi-dimensional.  We now give the proof of Prop.~\ref{inequality}.\par

\begin{proof}[Proof of Proposition~\ref{inequality}](1) We denote by $\mu:M\to \k^*\cong \k$ a moment map,
where $\k$ is the Lie algebra of $K$, identified with its dual
$\k^*$ by means of a $Ad(K)$-invariant scalar product
$<\cdot,\cdot>$ on $\k$. We fix a point $p\in M^L$ and, since
principal isotropy subgroups are conjugate, we may suppose that
$p$ belongs also to the open dense subset $M_\mu:=\{q\in M|\ \dim
K\cdot\mu(q)\geq\dim K\cdot\mu(w),\ \forall w\in M\}$. Then it is
known (see\ \cite{HW}, p.267) that, given a $K$-regular point
$p\in M$ with $K_p = L$, $L$ is normal in $K_{\mu(p)}$ with $K_{\mu(p)}/L$
abelian; moreover the kernel $\ker(d\mu_p)\cap T_p(K\cdot
p)$ is a $L$-invariant subspace of $T_p(K\cdot p)$, which can be
identified with the tangent space $T_p(K_{\mu(p)}\cdot p)$.
This implies that
$$\dim (\ker(d\mu_p)\cap T_p(K\cdot p)) = \rk(K) - \rk(L)$$
and also that $\ker(d\mu_p)\cap T_p(K\cdot p)\subset (T_pM)^L$.
Since the normal space $\nu_p = (T_p(K\cdot p))^\perp$ with
respect to any $K$-invariant scalar product on $M$ is obviously
contained in $(T_pM)^L$ because $p$ is regular, we obtain our first claim,
recalling that $\dim T_pM^L=\dim \overline M.$\\
(2) The fixed point set $(K\cdot p)^L$ clearly identifies with the
orbit $N_K(L)\cdot p$ and the equality  holds if and only if
$\dim(T_p(N_K(L)\cdot p)) = \dim K_{\mu(p)}\cdot p$, i.e. if and only if
$(N_K(L)/L)^o = K_{\mu(p)}/L$, which is abelian.\\
(3) Suppose that $\rk(K) =
\rk(L)$; this implies that $L=K_{\mu(p)}$ because $K_{\mu(p)}/L$
is abelian. Since $L$ has maximal rank, $V^L\cap T_p(K\cdot p) =
\{0\}$; moreover since $V^L$ contains the normal space of the
orbit $K\cdot p$ at $p$, we see that $V^L$ coincides with
$(T_pK\cdot p)^\perp$ and therefore the orbit $K\cdot p$ is
complex. In particular, since the $K$-action is Hamiltonian,
the group $K$ is semisimple by the same arguments used in the
proof of Prop.~\ref{flag}.\\
If $M=V$ is a linear space, then $K\cdot p$ is a compact
connected complex submanifold, hence a point and the $K$-action is trivial.\\
Suppose now that $M$ is a locally irreducible K\"ahler manifold.
We write $\k = \l + \m$, an $\Ad(L)$-invariant
decomposition, where $\m$ identifies with $T_p(K\cdot p)$; we
decompose $\m = \bigoplus_i \m_i$ into irreducible
$\Ad(L)$-invariant subspaces, which are mutually inequivalent,
because $K/L$ is a generalized flag manifold of a semisimple
compact Lie group (see e.g.~\cite{W}).
If $\xi\in V^L$, then the shape operator $A_\xi$, viewed as an
element of $\operatorname{End}(\m)$, commutes with $\Ad(L)$ and
therefore maps each $\m_i$ into itself; being self-adjoint, we
have that $A_\xi|_{\m_i} = \lambda_i I$ for some $\lambda_i\in\R$
by Schur's Lemma. Moreover, the complex structure $J$ of $V$
leaves each $\m_i$ invariant and therefore commutes with $A_\xi$;
on the other hand it is well-known (see~\cite{KN}) that $A_\xi$
anti-commutes with $J$, hence $A_\xi = 0$ for all $\xi\in V^L$,
meaning that every principal orbit is totally geodesic. The
open subset of principal points admits two complementary, totally
geodesic foliations, given by the $K$-orbits and by the fixed point sets of
principal isotropies; it is well-known that two complementary
totally geodesic foliations are parallel, contradicting the
local irreducibility of $M$.
\end{proof}
\begin{defi} A proper $K$-Hamiltonian action on a K\"ahler manifold $M$ by holomorphic isometries is called
complex asystatic (or $\C$-asystatic) if $$\dim \overline M=c+\rk
K-\rk L$$ or, equivalently, if $(N_K(L)/L)^o$ is abelian.
\end{defi}
\noindent{\bf Example.} On the quadric $Q_n =
\SO(n+2)/\SO(2)\times \SO(n)$, we consider the action of
$K=\SO(2)\times \SO(n)$. A principal isotropy subgroup $L$ is
given by ${\mathbb Z}_2\times \SO(n-2)$ and $(N_K(L)/L)^o =
(\SO(2))^2$, hence the action is $\C$-asystatic.
\medskip\\
\begin{rem} The example $K=\SO(n)$ acting on $\C^n$ shows that $N_K(L)/L$
is in general not abelian, even if the action is $\C$-asystatic;
indeed in this case $N_K(L)/L=T^1\cdot \Z_2$.
\end{rem}\\\\
The following proposition shows that the condition of $\C$-asystatic is
preserved in the slice representation at complex orbits.
\begin{prop} If the
$K$-action on $M$ is $\C$-asystatic then the slice representation
at a complex $K$-orbit is $\C$-asystatic.
\end{prop}
\begin{pf} Given a complex orbit $K\cdot x$, we have $\rk(K)=\rk(K_x)$; let $L\subset K_x$ be
a principal isotropy subgroup and denote by $\nu_x$ the normal space at $x$ of the orbit $K\cdot x$.
$$\dim (\nu_x (T_x K\cdot x))^{K_p}\leq \dim\overline M=c+\rk K-\rk
K_p$$this is equal to $c+\rk (K_x)-\rk {(K_x)}_p$ and we get the
claim.
\end{pf}
We now close this section  with the following lemma that
might be useful in order to classify $\C$-asystatic actions.
\begin{lem}\label{extension} Let $K$ be a compact connected Lie group acting on a K\"ahler manifold $M$
and $K'\subset K$ a compact connected subgroup whose action on $M$
is $\C$-asystatic, then the action of $K$ is $\C$-asystatic,
provided one of the following conditions holds:
\begin{enumerate}
\item $K$ and $K'$ have the same orbits;
\item $K = K'\cdot \operatorname{T}^1$, where $\operatorname{T}^1$ is a one dimensional torus.
\end{enumerate}
\end{lem}
\begin{pf} If (1) holds, let $N:= K/L = K'/L'$ be a common principal orbit through a point $p$ with $L'\subset L$; if $c$ denotes the common cohomogeneity
then
$$c + \rk(K) - \rk(L) \leq \dim (T_p M)^{L} \leq (T_p M)^{L'} = c + \rk(K') - \rk(L').$$
Our claim follows from the fact that $\rk(K) - \rk(L) = \rk(K') - \rk(L')$ (see e.g.\ \cite{O}, p. 207).\\
If (2) holds, we consider the moment maps $\mu,\mu'$ relative to
$K,K'$ respectively; the Lie algebra $\k$ of $K$ splits as $\k =
\k'\oplus \R\cdot Z$, where $Z$ is a generator of the Lie algebra
of $\operatorname{T}^1$ and the moment map $\mu'=\pi\circ\mu$,
where $\pi:\k\to\k'$ is the projection. Let $p\in M$ be a
principal point, then $\mu(K\cdot p)=K/H$, where $H\subset K$ is
the centralizer of a torus, hence of the form
$H=H'\cdot\operatorname{T}^1$ with $H'\subset K'$ the centralizer
of some torus in $K'$; moreover $\mu'(K'\cdot
p)=\pi\circ\mu(K'\cdot p)=\pi(K/H)=K'/H'$, because $K'$ acts
transitively on $K/H$.\\ Let $L'\subseteq L$ be principal isotropy
subgroups of $K',K$ respectively; by (1), we can suppose that
$K'\cdot p \neq K\cdot p$, hence $\l = \l'$, where $\l,\l'$ are
the Lie algebras of $L,L'$ resp. We decompose $\k=\l+\a+\n$, where
$\l+\a=\h$, $Z\in \a$ and $[\l,\a]=0$ and $\n=(\l+\a)^\perp$
w.r.t. an $\Ad(K)$-invariant scalar product $\langle
\cdot,\cdot\rangle$ with $\langle\k',Z\rangle = 0$. Now,
$\n^{\Ad(L)}\subseteq \n^{\Ad(L')} =\{0\}$ because the $K'$-action
is $\C$-asystatic; it then follows that
$$\dim (T_pK\cdot p)^L = \dim \a = \dim (\h'/\l) + 1 =$$$$= \rk(K')-\rk(L)+1 = \rk(K)-\rk(L)$$
hence our claim. \ \end{pf}

\begin{rem} We remark that if $K'\subset K$ have the same orbits and the $K$-action is $\C$-asystatic, then the
${K'}$-action is not necessarily $\C$-asystatic, as the example of
the linear action of $K'=\Sp(n)\subset K=\UU(2n)$ shows.\end{rem}

\section{$\C$-asystatic actions and spherical orbits}
In this section we prove our main Theorem \ref{spherical}. We
recall here that a complex homogeneous space $X = G/H$, where $G$
is a complex algebraic group, is called {\em spherical} if a Borel
subgroup of $G$ has an open orbit in $X$.  We recall also that  if
$(M,g)$ is a K\"ahler manifold with K\"ahler form $\omega$ and $K$
is a compact connected Lie subgroup of its full isometry group,
then the $K$-action is called {\em coisotropic} if the principal
$K$-orbits are coisotropic with respect to $\omega$.\par
\begin{proof}[Proof of Theorem~\ref{spherical}] Let $K$ be a maximal compact Lie subgroup
of $G$ whose action on  $M$ is Hamiltonian.
 Let $\mu:M\to \k^*$ be a moment map for the $K$-action and let
$M_\mu := \{v\in M;\ \dim K(\mu(v)) \geq \dim K(\mu(w))\ \forall\
w\in M\}$. Since $M_\mu$ is an open set, the intersection
$M_\mu\cap M_{\operatorname{reg}}$ is non empty and open; we fix
$p$ in $M_\mu\cap M_{\operatorname{reg}}$ and we denote by
${\mathcal O}$ the orbit $G\cdot p$ under the complexified group
$G$. We now prove that the orbit $K\cdot p$ is coisotropic in
$\mathcal O$, hence every principal orbit of $K$ in $\mathcal O$
will be coisotropic by Theorem $3$ in~\cite{HW}. By a result due
to Wurzbacher (Theorem 2.1, p. 542~\cite{Wu}) this means that
the $G$-orbit ${\mathcal O}$ is spherical.\\
In order to prove our claim, we denote by $\a_p$ the tangent space
$T_p(K_{\mu(p)}\cdot p)$ and by $\n_p$ the orthogonal complement $\n_p
= \a_p^\perp \cap T_p(K\cdot p)$; by the assumption of
$\C$-asystatic, the fixed point set $(T_p M)^L$ is given by $(T_p
M)^L = \nu_p \oplus\a_p$, where $\nu_p$ denotes the normal space $(T_p(K\cdot p)^\perp$;
since $T_p M^L$ is a complex subspace,
we have that $\n_p$ is complex too. This means that $T_p\mathcal O$
is given by $T_p(K\cdot p) \oplus J(\a_p)$, where $J$ denotes the
complex structure of $M$, so that the normal space of $T_p(K\cdot
p)$ inside $T_p\mathcal O$ is $J\a_p$; hence the orbit $K\cdot p$ is
coisotropic inside $\mathcal O$.\\
This argument shows that for any point $x$ in the open set $\in M
_\mu\cap M_{\operatorname{reg}}$ the orbit $G\cdot x$ is
spherical; our claim now follows from~\cite{A1}, where it is
proved that every $G$-orbit is spherical if this happens in an
open dense subset.\
\end{proof}

In the following we will denote by $G/H$  a $G$-spherical orbit
and  by $N_G(H)$ the normalizer of $H$ in $G$, using~\cite{Vi}
(Proposition 10, p. 19), under the same hypotheses of Theorem
\ref{spherical}, we get
\begin{corollary} If the action of a compact Lie subgroup of $G$  on  $M$
is $\C$-asystatic then $N_G(H)/H$ is abelian.
\end{corollary}

\section{The linear Case}

All representations of a reductive algebraic group $G$ with
spherical orbits have been classified in~\cite{A2}. In order to
have the classification of $\C$-asystatic representations, we have
to go through the list in~\cite{A2} and select which are
$\C$-asystatic. Since the computation of a principal isotropy is
not easily accessible in the literature, we provide a table of all
irreducible representations with spherical orbits, where we also
indicate a principal isotropy $L$. We now give the proof of the
classification Theorem.
\begin{proof}[Proof of Theorem~\ref{class}] By Theorem ~\ref{spherical}, we know that the $G$-orbits are spherical.
By Arzhantsev's result~\cite{A2},
an irreducible representation of a reductive algebraic group has
spherical orbits if and only if it appears in Table I (p.~291 in
~\cite{A2}) or is obtained from any of these by a torus extension;
here in Table 1, we indicate a compact real form $K$ with the
corresponding representation, a principal isotropy subgroup $L$,
the difference $f=\rk(K)-\rk(L)$ and the dimension $d$ of the
fixed point set for the isotropy action of $L$ on the tangent
space of the orbit $K/L$. Clearly, the representation is
$\C$-asystatic if and only if $f=d$. According to
Lemma~\ref{extension}, we then need only consider the cases which
are not $\C$-asystatic and compute the same invariants $f$ and $d$
for $K\cdot \T^1$; this is done in Table 2 and the full
classification is then obtained by taking all representations
which appear in Table 1 or 2 and are $\C$-asystatic.\end{proof}
\begin{table}[h]
\begin{tabular}{|c|c|c|c|c|c|}
\hline n.&$K$ & $\rho$ & $L$ & $f$& $d$\\
\hline $1$ &$\SU(n)$ & $\rho_1$& $\SU(n-1)$ & $1$& $1$\\
\hline $2$ &$\SU(2n)$ & $\Lambda^2\rho_1$& $\SU(2)^n$ & $n-1$& $n-1$\\
\hline $3$ &$\SU(2n+1)$ & $\Lambda^2\rho_1$& $\SU(2)^n$ & $n$& $n$\\
\hline $4$ &$\SU(n),n\geq 3$ & $S^2\rho_1$& $(\Z_2)^{n-1}$ & $n-1$& $n-1$\\
\hline $5$ &$\SU(2)$ & $S^2\rho_1$& $\Z_2$ & $1$& $3$\\
\hline $6$ &$\SO(n)$ & $\rho_1$& $\SO(n-2)$ & $1$& $1$\\
\hline $7$ &$\Sp(n)$ & $\rho_1$& $\Sp(n-1)$ & $1$& $3$\\
\hline $8$ &$\SU(n)\times\SU(m), n<m$ & $\rho_1\otimes\rho_1$& $\UU(1)^{n-1}\times \SU(m-n)$ & $n$& $n$\\
\hline $9$ &$\SU(n)\times\SU(n)$ & $\rho_1\otimes\rho_1$& $\UU(1)^{n-1}$ & $n-1$& $n-1$\\
\hline $10$ &$\SU(2)\times\Sp(n), n\geq 2$ & $\rho_1\otimes\rho_1$& $\UU(1)\times\Sp(n-2)$ & $2$& $4$\\
\hline $11$ &$\SU(3)\times\Sp(n), n\geq 3$ & $\rho_1\otimes\rho_1$& $\Sp(n-3)$ & $4$& $29$\\
\hline $12$ &$\SU(3)\times\Sp(2)$ & $\rho_1\otimes\rho_1$& $\{e\}$ & $4$& $18$\\
\hline $13$ &$\SU(4)\times\Sp(2)$ & $\rho_1\otimes\rho_1$& $\Z_2$ & $5$& $25$\\
\hline $14$ &$\SU(n)\times\Sp(2), n\geq 5$ & $\rho_1\otimes\rho_1$& $\SU(n-4)\times\Z_2$ & $6$& $26$\\
\hline $15$ &$\Spin(7)$ & $\rho_3$& $\SU(3)$ & $1$& $1$\\
\hline $16$ &$\Spin(9)$ & $\rho_4$& $\SU(3)$ & $2$& $4$\\
\hline $17$ &$\Spin(10)$ & $\rho_4$& $\SU(4)$ & $2$& $2$\\
\hline $18$ &$\GG_2$ & $\rho_1$& $\SU(2)$ & $1$& $3$\\
\hline $19$ &$\EE_6$ & $\rho_1$& $\Spin(8)$ & $2$& $2$\\
\hline
\end{tabular}\vspace{1em}

\caption{Semisimple $K$ such that $G$ has spherical orbits.}
\end{table}
\medskip

\begin{table}[h]
\begin{tabular}{|c|c|c|c|c|c|}

\hline n.&$K$ & $\rho$ & $L$ & $f$& $d$\\
\hline $1$ &$\UU(2)$ & $S^2\rho_1$& $(\Z_2)^2$ & $2$& $2$\\
\hline $2$ &$\Sp(n)\times\T^1$ & $\rho_1\otimes\epsilon$& $\Sp(n-1)\times\T^1$ & $1$& $1$\\
\hline $3$ &$\UU(2)\times\Sp(n)$ & $\rho_1\otimes\rho_1$& $\T^1\times\Sp(n-2)$ & $3$& $5$\\
\hline $4$ &$\UU(3)\times\Sp(n), n\geq 3$ & $\rho_1\otimes\rho_1$& $\Z_2\times\Sp(n-3)$ & $6$& $30$\\
\hline $5$ &$\UU(3)\times\Sp(2)$ & $\rho_1\otimes\rho_1$& $\Z_2$ & $5$& $19$\\
\hline $6$ &$\UU(n)\times\Sp(2), n\geq 4$ & $\rho_1\otimes\rho_1$& $\UU(n-4)\times\Z_2$ & $6$& $26$\\
\hline $7$ &$\Spin(9)\times\T^1$ & $\rho_4\otimes\epsilon$& $\SU(3)\times\Z_2$ & $3$& $5$\\
\hline $8$ &$\GG_2\times\T^1$ & $\rho_1\otimes\epsilon$& $\SU(2)\cdot\Z_2$ & $2$& $2$\\
\hline
\end{tabular}\vspace{1em}
\caption{Non-semisimple $K$ such that $G_{ss}$ acts with spherical
orbits, but $K_{ss}$ acts non $\C$-asystatically}\label{Tten}
\end{table}

\newpage
\noindent{\bf Notation:} In the above Tables $\rho_i$ denote the
standard representations.

\end{document}